\documentclass{amsart}
\usepackage{mathrsfs}
\usepackage{amssymb}
\usepackage{latexsym}
\usepackage{yfonts}
\usepackage{natbib}

\usepackage{graphicx,color}

\usepackage[plainpages=false,colorlinks,hyperindex,bookmarksopen,linkcolor=red,citecolor=blue,urlcolor=blue]{hyperref}

\bibpunct{[}{]}{;}{n}{,}{,}

\newtheorem{te}{Theorem}

\newtheorem{os}{Remark}

\numberwithin{equation}{section}

\allowdisplaybreaks

\begin{document}

\title[Fractional operators]{Continuous random walks and\\ fractional powers of operators}

\author{Mirko D'Ovidio} 
\address{Dipartimento di Scienze di Base e Applicate per l'Ingegneria, Sapienza University of Rome, A. Scarpa 10 - 00161, Rome, Italy}
\email{mirko.dovidio@uniroma1.it}

\keywords{Fractional Laplacian, Stable process, Compound Poisson process.}

\date{\today}

\subjclass[2010]{60J35, 60J50, 60J75}

\begin{abstract}
We derive a probabilistic representation for the Fourier symbols of the generators of some stable processes.
\end{abstract}

\maketitle

\section{Introduction and main results}

The connection between fractional operator in space and diffusion with long jumps has been pointed out by many researchers (see for example  \cite{beghinMacci,meerSilbook,valdinoci09} and the references therein). It is well known that the compound Poisson process is a continuous time stochastic process with jumps which arrive, according to a Poisson process,  with specific probability law for the size. Our aim is to characterize the jumps distribution in order to obtain singular limit measure characterizing fractional powers of operators.  

Let $N(t)$, $t>0$ be a Poisson process with rate $\lambda>0$. Let $Y_j$, $0\leq j \leq n$ be $n+1$ i.i.d. random jumps such that $Y_j \sim Y$ for all $j$, where the symbol ''$\sim$'' stands for equality in law. It is well known that
\begin{equation}
Z_t = \sum_{j=0}^{N(t)} Y_j - \lambda t \mathbb{E}Y, \quad t>0
\end{equation}
is the compensated Poisson process with generator
\begin{equation}
(\mathcal{A}f)(x) = \lambda \int_{\mathbb{R}} \left( f(x+y) - f(x) - y f^\prime(x) \right) \nu_Y(dy) \label{generator-CP}
\end{equation}
where $\nu_Y: \Omega \subseteq \mathbb{R} \mapsto [0,1]$ is the density law of $Y \in \Omega$. The latter is quite familiar in the representation of the fractional power of the Laplacian. Indeed, the fractional Laplace operator can be defined pointwise:
\begin{equation}
-(-\triangle)^\alpha f(\mathbf{x}) =  \int_{\mathbb{R}^d} \Big( f(\mathbf{x}+\mathbf{y}) - f(\mathbf{x}) - \mathbf{y} \cdot \nabla f(\mathbf{x}) \mathbf{1}_{(|\mathbf{y}|\leq 1)} \Big)\frac{C_d(\alpha)\, \mathbf{dy}}{|\mathbf{y}|^{2\alpha +d}}
\end{equation}
where $C_d(\alpha)$ is a constant depending on $d$ and $\alpha \in (0,1)$, $f$ is a suitable test function, $C^2$ function with bounded second derivative for instance.

In this short paper, we construct continuous random walks with exponential and Gaussian jumps driven by pseudo-differential operators with Fourier multiplier $\Phi_\gamma(\boldsymbol{\xi})$ which converges to $|\boldsymbol{\xi} |^{\beta}$ with $\beta \in (0,2)$ as $\gamma \to 0$. In particular, we first consider the random jump $Y=\gamma e^X \in [\gamma, \infty)$ where $X \sim Exp(\alpha)$ and $\alpha,\gamma >0$. We have that $P\{Y \in A\} = \int_A \nu_Y(dy)$ with $\nu_Y(y) = \alpha \gamma^\alpha  y^{-\alpha -1}\mathbf{1}_{(y\geq \gamma)}$. By ''symmetrizing'', we get that
\begin{equation}
 \nu^*_Y(y) =  q\, \nu_Y(-y)\mathbf{1}_{(y \leq -\gamma)} + p\,\nu_Y(y)\mathbf{1}_{(y \geq \gamma)}\label{sym-law-no-gamma}
\end{equation} 
is the density of $Y^* = \epsilon\, Y$ with Rademacher law
$$P\{\epsilon = +1\}= p, \quad P\{\epsilon = -1\}=q$$
for the random variable $\epsilon$ where, we obviously assume that $p+q=1$. For $p=q$, formula \eqref{sym-law-no-gamma} takes the form
\begin{equation}
\nu^*_Y(y) = \frac{1}{2}\nu_Y(|y|) = \frac{\alpha \gamma^{\alpha}}{2} |y|^{-\alpha -1} \mathbf{1}_{(|y|\geq \gamma)}
\end{equation}
and $Y^*$ is therefore written as
\begin{equation}
Y^* = \left\lbrace \begin{array}{ll}
\gamma e^X, & \textrm{with probability }\; 1/2,\\
-\gamma e^X, & \textrm{with probability }\; 1/2.
\end{array} \right .
\end{equation}
We write the corresponding compound Poisson process as follows
\begin{equation}
A(t) = \sum_{j=0}^{N(t)} Y^*_j = \sum_{j=0}^{N(t)}  \epsilon_j\, Y_j = \sum_{j=0}^{N(t)}  \epsilon_j\, \gamma_j\, e^{X_j}, \quad t>0 \label{A-proc}
\end{equation}
with $X_j \sim X$, $\epsilon_j \sim \epsilon$ and $\gamma_j=\gamma$ for all $j=0,1,2,\ldots $. We also assume that all the random variables we are dealing with are taken to be independent, that is $\mathbb{E}\epsilon \epsilon_j = \mathbb{E}\epsilon_j \epsilon_{j^\prime} = 0$, $\forall\, j,j^\prime$ such that $j\neq j^\prime$. Observe that, for all $\varepsilon >0$, $P\{ \gamma e^X < \varepsilon \} \to P\{ e^X < +\infty \} =1$ as $\gamma \to 0$. 

We recall that the symbol $\xrightarrow{d}$ stands for ''converges in distribution'' and state the following results.
\begin{te}
\label{zero-the}
Let $\mathfrak{H}^\alpha_j(t)$, $t>0$, $j=1,2$ be two independent stable subordinators. 
For given $p,q \geq 0$ such that $p+q=1$, $\alpha \in (0,1)$,
\begin{equation}
A(t/\gamma^{\alpha}) \xrightarrow[\gamma \to 0]{d} \mathfrak{H}^\alpha_{1}(pt^*) - \mathfrak{H}^\alpha_{2}(qt^*) \label{A-zero-thm}
\end{equation}
with $t^* : t \mapsto \lambda \Gamma(1-\alpha) t$ and generator
\begin{equation}
\mathcal{A}f(x) = - \lambda \Gamma(1-\alpha) \left( p  \frac{d^\alpha}{d x^\alpha}  + q \frac{d^\alpha}{d (-x)^\alpha} \right) f(x). \label{op-zero-thm}
\end{equation}
\end{te}
The Weyl's fractional derivatives appearing in \eqref{op-zero-thm} are defined as follows:
\begin{align*}
\frac{d^\alpha f}{d x^\alpha}(x) = & \frac{1}{\Gamma(1-\alpha)} \frac{d}{dx} \int_0^\infty f(x-y) \frac{dy}{y^{\alpha}}\\ 
= & \frac{\alpha}{\Gamma(1 - \alpha)}\int_0^\infty \left( f(x) - f(x-y) \right)\frac{dy}{y^{\alpha +1}};
\end{align*}
\begin{align*}
\frac{d^\alpha f}{d (-x)^\alpha}(x) = & \frac{-1}{\Gamma(1-\alpha)} \frac{d}{dx} \int_0^\infty f(x+y) \frac{dy}{y^{\alpha}}\\ 
= & \frac{\alpha}{\Gamma(1 - \alpha)}\int_0^\infty \left( f(x) - f(x+y) \right)\frac{dy}{y^{\alpha +1}}.
\end{align*}
We consider ''good'' functions $f : \mathbb{R} \mapsto [0,1]$ whereas, we obtain the Riemann-Liouville derivatives (left-handed and right-handed respectively) by considering $f: \mathbb{R}_{+} \mapsto [0,1]$ and $f: \mathbb{R}_{-} \mapsto [0,1]$ respectively. In the integrals above, for example, one can write $f(z)\mathbf{1}_{(z\geq 0)}$ and  $f (z) \mathbf{1}_{(z \leq 0)}$ and obtain the operator governing the (totally) positively and negatively skewed stable process, that is $q=0$ and $p=0$ respectively in \eqref{op-zero-thm}.

\begin{os}
A stable subordinator is a one-dimensional non-decreasing L\'{e}vy process with $\mathbb{E}\exp(-\mu \mathfrak{H}^\alpha_t) = \exp(-t \mu^\alpha)$, $\alpha \in (0,1)$. We observe that $\mathfrak{H}^\alpha_j(0)=0$ for $j=1,2$ and therefore the process \eqref{A-zero-thm} converges (in distribution) to a totally positively (if $q=0$, $P\{\epsilon_j =+1\}=1$, $\forall\, j$) or negatively (if $p=0$, $P\{\epsilon_j =-1\}=1$, $\forall\, j$) skewed stable process. Furthermore, for a given process $X_t$, $t>0$ we notice that $X_{\theta t}$ runs slower than $X_t$ as well as $\theta$ is less that $1$.
\end{os}

\begin{te}
\label{first-the}
Let $\mathfrak{S}^{\beta}(t)$, $t>0$ be a symmetric stable process with $\beta \in (0,2)$. Then, for $p=q=1/2$, $\alpha \in (0,2)$,  
\begin{equation}
A(t/\gamma^{\alpha}) \xrightarrow[\gamma \to 0]{d} \mathfrak{S}^{\alpha}(t^*) \label{A-first-thm}
\end{equation}
with $t^* : t \mapsto \alpha \lambda C t$ and generator
\begin{equation}
\mathcal{A}f(x) = - \alpha \lambda C \frac{d^\alpha f}{d |x|^\alpha}(x)
\end{equation}
where
\begin{equation}
C = \frac{1}{2} \int_\mathbb{R} \frac{1-\cos y}{|y|^{\alpha +1}}dy.
\end{equation}
\end{te}

We now introduce the reciprocal gamma random variable $E_\alpha$, $\alpha>0$  with $P\{E_\alpha \in dx\} = \frac{x^{-\alpha -1}}{\Gamma(\alpha)} e^{-1/x}dx$ (the reciprocal gamma process has interesting connections with stable subordinators and Bessel processes, see for example \cite{Dov3}). We also consider the (normal) random vector $\mathbf{Y} \sim N(\mathbf{0}, \sigma^2_\alpha)$, $\mathbf{Y} \in \mathbb{R}^d$, with random variance $\sigma^2_\alpha \sim \frac{\gamma}{2}E_\alpha$ for some $\gamma >0$ and define the process
\begin{equation}
\mathbf{A}(t) =  \sum_{j=0}^{N(t)} \epsilon_j \mathbf{Y}_j
\end{equation} 
where $\epsilon_j \sim \epsilon$ $\forall\, j$, $\epsilon$ has Rademacher law as above, $\mathbf{Y}_j \sim \mathbf{Y}$ $\forall\, j$ and $\mathbb{E}\epsilon_j \epsilon_{j^\prime} = 0$ for all $j, j^\prime$ such that $j \neq j^\prime$. We notice that 
\begin{align*}
P\{ N(0, \sigma^2_\alpha) > x \} \leq   \frac{1}{x} \int_x^\infty y  \left(\int_0^\infty \frac{e^{-\frac{y^2}{2s}}}{\sqrt{2 s}} P\{ \sigma^2_\alpha \in ds \} \right) dy  \approx \gamma^\alpha x^{-2\alpha}
\end{align*}
and therefore, for large $x$,
\begin{equation*}
P\{N(0, \sigma^2_\alpha) \in dx\}/dx \approx 2\alpha \gamma^\alpha x^{-2\alpha -1}.
\end{equation*}
After some calculations we explicitly write the law of $\mathbf{Y} \in \mathbb{R}^d$ as follows
\begin{equation}
\nu_{\mathbf{Y}}(\mathbf{y}) = \frac{\Gamma(\alpha + \frac{d}{2})}{\pi^\frac{d}{2}\, \Gamma(\alpha) } \frac{\gamma^\alpha}{(|\mathbf{y}|^2 + \gamma)^{\alpha + \frac{d}{2}}}, \quad \mathbf{y} \in \mathbb{R}^d. \label{nu-Y-bold}
\end{equation}

We are now ready to present the next result.
\begin{te}
\label{second-the}
Let $\boldsymbol{\mathfrak{S}}^{\beta}(t) \in \mathbb{R}^d$, $t>0$ be an isotropic stable process with $\beta \in (0,2)$. Then, for $\alpha \in (0,1)$,  
\begin{equation}
\mathbf{A}(t/\gamma^{\alpha}) \xrightarrow[\gamma \to 0]{d} \boldsymbol{\mathfrak{S}}^{2\alpha}(t^*) \label{A-second-thm}
\end{equation}
with $t^* : t \mapsto \lambda \frac{\Gamma(\alpha + \frac{d}{2})}{\pi^\frac{d}{2}\, \Gamma(\alpha)} C t$ and infinitesimal generator
\begin{equation}
\mathcal{A}f(\mathbf{x}) = - \lambda \frac{\Gamma(\alpha + \frac{d}{2})}{\pi^\frac{d}{2}\, \Gamma(\alpha)} C (-\triangle)^\alpha f(\mathbf{x})
\end{equation}
where
\begin{equation}
C =  \int_{\mathbb{R}^d} \frac{1- \cos y_1}{|\mathbf{y}|^{2\alpha + d}} \mathbf{dy}.
\end{equation}
\end{te}

\section{Compensated Poisson and fractional Laplace operator}
Let us consider the L\'{e}vy process $\mathbf{F}_t$, $t>0$, with associated Feller semigroup $T_t \, f(\mathbf{x}) = \mathbb{E}f(\mathbf{F}_t - \mathbf{x})$ solving $\partial_t u = \mathcal{A}u$, $u_0=f$. The operator $\mathcal{A}$ is the infinitesimal generator of $\mathbf{F}_t$, $t>0$ and the following representation holds
\begin{equation}
(\mathcal{A}f)(\mathbf{x}) = -\frac{1}{(2\pi)^d} \int_{\mathbb{R}^d} e^{-i \boldsymbol{\xi}\cdot \mathbf{x}} \Phi(\boldsymbol{\xi}) \widehat{f}(\boldsymbol{\xi})\boldsymbol{d \xi} \label{inv-fourier-symbol}
\end{equation}
for all functions in the domain 
\begin{equation}
D(\mathcal{A}) = \left\lbrace f \in L^2(\mathbb{R}^d,\mathbf{dx})\,:\, \int_{\mathbb{R}^d} \Phi(\boldsymbol{\xi}) |\widehat{f}(\boldsymbol{\xi})|^2 \boldsymbol{d \xi}< \infty \right\rbrace
\end{equation}
where $\widehat{f}(\boldsymbol{\xi}) = \int_{\mathbb{R}^d} e^{i \boldsymbol{\xi}\cdot \mathbf{x}} f(\mathbf{x}) \mathbf{dx}$ is the Fourier transform of $f$, $\Phi(\cdot)$ is continuous and negative definite. We say that $T_t$ is a pseudo-differential operator with symbol $\exp(-t\Phi)$ and, $\Phi$ is the Fourier multiplier (or {\bf Fourier symbol}) of $\mathcal{A}$, $\widehat{(\mathcal{A}f)}(\boldsymbol{\xi}) = -\Phi(\boldsymbol{\xi}) \widehat{f}(\boldsymbol{\xi})$. Furthermore (as in \cite{Jacob1998}), we write
\begin{equation}
- \partial_t\, \mathbb{E}e^{i\boldsymbol{\xi}\cdot \mathbf{F}_t} \Big|_{t=0} = \Phi(\boldsymbol{\xi}). 
\end{equation}
It is well known that, for $\Phi(\boldsymbol{\xi})=|\boldsymbol{\xi} |^\alpha$, formula \eqref{inv-fourier-symbol} gives us the fractional power of the Laplace operator which can be also expressed as 
\begin{align}
-(-\triangle)^{\alpha} f(\mathbf{x}) = & C_d(\alpha) \, \textrm{p.v.}  \int_{\mathbb{R}^d} \frac{f(\mathbf{y}) -f(\mathbf{x})}{|\mathbf{x}-\mathbf{y}|^{2\alpha +d}}\mathbf{dy} \notag \\
= & C_d(\alpha) \, \textrm{p.v.} \int_{\mathbb{R}} \frac{f(\mathbf{x}+\mathbf{y}) -f(\mathbf{x})}{|\mathbf{y}|^{2\alpha +d}}\mathbf{dy} \label{lap-princ-val}
\end{align}
where ''p.v.'' stands for the ''principal value'' of the singular integrals above near the origin. For $\alpha \in (0,1)$, the fractional Laplace operator can be defined, for $f \in \mathscr{S}$ (the space of rapidely decaying $C^\infty$ functions), as follows
\begin{align}
-(-\triangle)^\alpha f(\mathbf{x}) = & \frac{C_d(\alpha)}{2} \int_{\mathbb{R}^d} \frac{f(\mathbf{x}+\mathbf{y}) + f(\mathbf{x}-\mathbf{y}) - 2f(\mathbf{x})}{|\mathbf{x}-\mathbf{y}|^{2\alpha +d}}\mathbf{dy}\notag\\
= & \frac{C_d(\alpha)}{2} \int_{\mathbb{R}^d} \frac{f(\mathbf{x}+\mathbf{y}) + f(\mathbf{x}-\mathbf{y}) - 2f(\mathbf{x})}{|\mathbf{y}|^{2\alpha +d}}\mathbf{dy}, \quad \forall\, \mathbf{x} \in \mathbb{R}^d. 
\end{align}
This representation comes out by considering straightforward calculations and removes the singularity at the origin (\cite{NPV2011}). Indeed, from the second order Taylor expansion of the smooth function $f$ ($f \in \mathscr{S}$) we obtain
\begin{equation}
\frac{f(\mathbf{x}+\mathbf{y}) + f(\mathbf{x}-\mathbf{y}) - 2f(\mathbf{x})}{|\mathbf{y}|^{2\alpha +d}} \leq \frac{\| D^2 f\|_{L^\infty}}{|\mathbf{y}|^{2\alpha + d - 2}}
\end{equation}
which is integrable near the origin provided that $\alpha \in (0,1)$. The constant $C_d(\alpha)$ must be considered in order to obtain $\widehat{(-\triangle)^\alpha f(\cdot) }(\boldsymbol{\xi}) = |\boldsymbol{\xi}|^\alpha \widehat{f}(\boldsymbol{\xi})$. 

\begin{os}
Let $\mathfrak{H}^\alpha_t$, $t>0$ be a stable subordinator. The generator of $\mathbf{F}_{\mathfrak{H}^\alpha_t}$, $t>0$ is given by the beautiful formula
\begin{equation}
-(-\mathcal{A})^\alpha f(x) = \frac{\alpha}{\Gamma(1-\alpha)} \int_0^\infty \Big( T_s \, f(x) - f(x) \Big) \frac{ds}{s^{\alpha +1}}
\end{equation}
for all $f \in \mathscr{S}$ ($T_s=e^{s\mathcal{A}}$ is the Feller semigroup of $\mathbf{F}_t$, $t>0$). 
\end{os}

Formula \eqref{generator-CP} can be obtained by considering the following characteristic function
\begin{align*}
\mathbb{E} e^{i\xi Z_t} = & \mathbb{E} \prod_{j=0}^{N(t)} e^{i\xi Y_j} \,  e^{-i\xi \alpha t \mathbb{E}Y} \\
= & \mathbb{E} \left(  \mathbb{E} e^{i\xi Y} \right)^{N(t)} \, e^{-i\xi \alpha t \mathbb{E}Y} \\
= & \exp\Big( \lambda t \mathbb{E} \left(e^{i\xi Y} - 1 -i \xi Y \right) \Big).
\end{align*}
Therefore, we get that 
\begin{equation*}
\partial_t\, \mathbb{E} e^{i\xi Z_t} \Big|_{t=0} = \lambda  \mathbb{E} (e^{i\xi Y} - 1 -i \xi Y) = \lambda \int_{\mathbb{R}} (e^{i\xi y} - 1 -i \xi y)\nu_Y(dy) = -\Phi(\xi).
\end{equation*}
If $Y_j \sim Y$ are symmetric random variables such that $\mathbb{E}Y_j=\mathbb{E}Y=0$ for all $j=1,2,\ldots$, than $\nu_Y(y) = \nu_Y(-y)$ and 
\begin{equation}
 \int_{\mathbb{R}} y f^\prime(x)\, \nu_Y(dy)  = f^\prime(x) \int_{\mathbb{R}\setminus B_r} y\, \nu_Y(dy) + f^\prime(x) \int_{B_r} y\, \nu_Y(dy)  = 0 \label{int-zero}
\end{equation}
where we also include those density law $\nu_Y(\cdot)$ for which \eqref{int-zero} holds as principal value. If \eqref{int-zero} holds true, then formula \eqref{generator-CP} takes the form
\begin{equation*}
(\mathcal{A}f)(x) = \lambda \int_{\mathbb{R}} \left( f(x+y)-f(x) \right) \nu_Y(dy)
\end{equation*}
and the integral converges depending on $\nu_Y(\cdot)$. If we choose $\nu_Y(dy) = 2\alpha |y|^{-2\alpha -1} dy$ for instance, then the integral must be understood in the principal value sense and we get the fractional Laplace operator as formula \eqref{lap-princ-val} entails.

\section{Proof of Theorem \ref{zero-the}}
The characteristic function of \eqref{A-proc} is written as follows
\begin{align*}
\mathbb{E} \exp\left( i\xi  \sum_{j=0}^{N(t)} \epsilon_j\, Y_j \right) = & \mathbb{E} \left( \mathbb{E} e^{i\xi \epsilon Y} \right)^{N(t)}\\
= & \exp \left( \lambda t (\mathbb{E} e^{i\xi \epsilon  Y} - 1) \right)\\
= & \exp \Bigg( \lambda t \Big(p \mathbb{E} e^{i\xi   Y} + q \mathbb{E} e^{-i\xi  Y} - (p+q) \Big) \Bigg)\\
= & \exp \Bigg( \lambda t \Big(p (\mathbb{E} e^{i\xi   Y} - 1) + q (\mathbb{E} e^{-i\xi  Y} - 1) \Big) \Bigg).
\end{align*}
From this, we immediately  get
\begin{equation*}
\mathbb{E} \exp\left( i\xi A(t/\gamma^{\alpha}) \right) = \exp \Bigg( \frac{\lambda t}{\gamma^\alpha} \Big(p (\mathbb{E} e^{i\xi   Y} - 1) + q (\mathbb{E} e^{-i\xi  Y} - 1) \Big) \Bigg).
\end{equation*}
We recall that the L\'{e}vy symbol of a stable subordinator is a mapping from $\mathbb{R} \mapsto \mathbb{C}$ which takes the form
\begin{equation}
-(-i\xi)^\alpha = \frac{\alpha}{\Gamma(1-\alpha)} \int_0^\infty \left(e^{i\xi y} - 1 \right) \frac{dy}{y^{\alpha+1}} \label{fur-symb-subordinator}
\end{equation}
for $\alpha \in (0,1)$. The Fourier symbol (depending on $\gamma$) of the characteristic function of $A(t/\gamma^{\alpha})$ is therefore given by
\begin{align*}
\Phi_\gamma(\xi) = & - \partial_t\, \mathbb{E} \exp\left( i\xi  A(t/\gamma^{\alpha}) \right) \Bigg|_{t=0}\\
= & - \frac{\lambda }{\gamma^\alpha} \Big(p \int_0^\infty (e^{i\xi y} - 1) \nu_Y(dy) + q \int_0^\infty (e^{-i\xi y} - 1)\nu_Y(dy) \Big) \\
= & -\lambda  p \int_\gamma^\infty (e^{i\xi y} - 1) \frac{\alpha dy}{y^{\alpha +1}} - \lambda q \int_\gamma (e^{-i\xi y} - 1)\frac{\alpha dy}{y^{\alpha +1}}  .
\end{align*}
For $\gamma \to 0$ we obtain
\begin{align}
\Phi_\gamma(\xi) \to \Phi(\xi)=  &  \lambda \Gamma(1-\alpha )  \Big(p (-i\xi)^\alpha + q (i\xi)^\alpha \Big), \quad \alpha \in (0,1)  \label{phi-symb-zero}
\end{align}
and, from \eqref{inv-fourier-symbol} we arrive at
\begin{equation}
\mathcal{A}f(x) = - \lambda \Gamma(1-\alpha) \left( p  \frac{d^\alpha}{d x^\alpha} + q \frac{d^\alpha}{d (-x)^\alpha}  \right) f(x).
\end{equation}
The fact that the L\'{e}vy process \eqref{A-zero-thm} has infinitesimal generator \eqref{op-zero-thm} comes directly from the characteristic function
\begin{align*}
\mathbb{E} \exp\Big(i\xi \mathfrak{H}^\alpha_{1}(pt^*) - i\xi \mathfrak{H}^\alpha_{2}(qt^*)  \Big) = \exp\left(-t^* p(-i\xi)^\alpha-t^*q(i\xi)^\alpha \right) 
\end{align*}
where $t^* = \lambda \Gamma(1-\alpha) t>0$. Thus, we get that
\begin{equation*}
-\partial_t \mathbb{E}\exp\Big(i\xi \mathfrak{H}^\alpha_{1}(pt^*) - i\xi \mathfrak{H}^\alpha_{2}(qt^*)  \Big) \Big|_{t=0} = \lambda \Gamma(1-\alpha )  \Big(p (-i\xi)^\alpha + q (i\xi)^\alpha \Big)
\end{equation*}
which coincides with $\Phi(\xi)$ in \eqref{phi-symb-zero}.

In the last calculations we have used the fact that
\begin{equation*}
\mathbb{E}e^{i\xi \mathfrak{H}^\alpha(t)} = \exp\left(-t (-i\xi)^\alpha \right) = \exp\left( -t |\xi|^{\alpha} e^{-i\frac{\pi \alpha}{2} \frac{\xi}{|\xi |}} \right), \quad \xi \in \mathbb{R}, \; t\geq 0
\end{equation*}
and thus,
\begin{equation*}
\mathbb{E}e^{-i\xi \mathfrak{H}^\alpha(t)} = \exp\left( -t |\xi|^{\alpha} e^{i\frac{\pi \alpha}{2} \frac{\xi}{|\xi |}} \right) = \exp\left( -t (i\xi)^\alpha \right).
\end{equation*}
The Fourier transforms of the Weyl's fractional derivatives, for $\alpha \in (0,1)$, are given by (\cite{SKM93})
\begin{equation}
\int_{\mathbb{R}} e^{i\xi x} \frac{d^\alpha}{d (\pm x)^{\alpha}} f(x)dx = (\mp i \xi)^\alpha \widehat{f}(\xi), \quad f \in L^1(\mathbb{R}). \label{fur-weyl-der}
\end{equation}

\section{Proof of Theorem \ref{first-the}}
For $p=q=1/2$ and $\alpha \in (0,2)$ we obtain that
\begin{align*}
\Phi_\gamma(\xi) = & -\partial_t\, \mathbb{E} \exp\left( i\xi  A(t/\gamma^{\alpha}) \right) \Bigg|_{t=0}\\
=  & -\frac{\lambda }{2 \gamma^{\alpha}} (\mathbb{E} e^{i\xi  Y} + \mathbb{E} e^{-i\xi Y} - 2)\\
= & -\frac{\lambda}{2 \gamma^{\alpha}} \int_{0}^{\infty} \left(e^{i\xi  y} +  e^{-i\xi y} - 2 \right) \nu_Y(dy)\\
= & -\frac{\lambda}{\gamma^{\alpha}} \int_{\mathbb{R}} \left( \cos(\xi y) - 1\right) \nu^*_Y(dy)
\end{align*}
where, we recall that
\begin{equation*}
\nu^*_Y(y) = \frac{\alpha \gamma^{\alpha}}{2} |y|^{-\alpha -1} \mathbf{1}_{(|y|\geq \gamma)}.
\end{equation*}
We explicitly have that
\begin{align*}
\Phi_\gamma(\xi) = & -\frac{\lambda}{\gamma^{\alpha}} \int_{\mathbb{R}} \left( \cos(\xi y) - 1 \right) \nu^*_Y(dy)\\
= & -\frac{\alpha \lambda}{2}  \int_{\mathbb{R}\setminus B_\gamma} \left( \cos(\xi y) - 1 \right) |y|^{-\alpha -1} dy.
\end{align*}
By taking the limit for $\gamma \to 0$, we obtain 
\begin{equation}
\Phi_\gamma(\xi) \to \Phi(\xi) = -\frac{\alpha \lambda}{2}  \int_{\mathbb{R}} \left( \cos(\xi y) - 1 \right) |y|^{-\alpha -1} dy =  \alpha \lambda C |\xi |^{\alpha} \label{phi-first-limit}
\end{equation}
where, due to the fact that $\left( \cos(y) - 1 \right) |y|^{-\alpha -1} \leq  y^2 |y|^{-\alpha-1}$ by Taylor expansion near the origin, we obtain that
\begin{equation*}
0 < C = \frac{1}{2}\int_\mathbb{R} \left( 1 - \cos(y)\right) |y|^{-\alpha -1} dy < \infty 
\end{equation*}
and $|\xi |^{\alpha}$ is the Fourier multiplier of the infinitesimal generator of a stable symmetric process. Indeed,  for the symmetric stable process $\mathfrak{S}^\beta(t)$, $t>0$, $\beta \in (0,2)$, we have that 
\begin{equation*}
-\partial_t\, \mathbb{E}e^{i\xi \mathfrak{S}^{\beta}(t)} \Big|_{t=0} =  |\xi |^{\beta}
\end{equation*}
and \eqref{inv-fourier-symbol} holds with 
\begin{equation}
\mathcal{A}f(x)= - \frac{d^\beta f}{d |x|^\beta}(x) = - \frac{\sigma}{2} \left( \frac{d^\beta f}{d x^\beta}(x) + \frac{d^\beta f}{d (-x)^\beta}(x) \right) \label{riesz-op}
\end{equation}
where $\sigma = (\cos \pi \beta /2)^{-1}$. The Fourier symbol of the Riesz operator \eqref{riesz-op} is written as (see formula \eqref{fur-weyl-der})
\begin{align*}
\int_{\mathbb{R}} e^{i\xi x} \frac{d^\beta f}{d |x|^\beta}(x) dx = \frac{\sigma}{2} \left( (-i \xi)^\beta + (i\xi)^\beta \right) \widehat{f}(\xi)= |\xi |^\beta\, \widehat{f}(\xi).
\end{align*}
Therefore, from \eqref{phi-first-limit}, we conclude that
\begin{equation*}
\sum_{j=0}^{N(t/\gamma^{\alpha})} Y^*_j  \stackrel{\gamma \to 0}{\longrightarrow} \mathfrak{S}^{\alpha}(t^*)
\end{equation*}
in distribution, where
\begin{equation*}
\mathbb{E}e^{i\xi \mathfrak{S}^{\alpha}(t^*)} = \exp\left( - t^*  |\xi |^{\alpha} \right)
\end{equation*}
and $t^* = \alpha \lambda C t$, $t>0$. Furthermore, the generator of  $\mathfrak{S}^\alpha(t^*)$ with $\alpha \in (0,2)$ is 
\begin{equation*}
-\alpha \lambda C \frac{d^\alpha f}{d |x|^\alpha}(x)= - \frac{1}{2\pi} \int_{\mathbb{R}} e^{-i\xi x} \Phi(\xi) \widehat{f}(\xi)d\xi.
\end{equation*}

We notice that, for $\alpha \in (0,1)$, the constant $\alpha \lambda C$ equals
\begin{align*}
\frac{\alpha \lambda}{2}\int_\mathbb{R} \left( 1 - \cos(y)\right) |y|^{-\alpha -1} dy = & 2 \int_0^\infty  \left( 1 - \cos(y)\right) \frac{dy}{y^{\alpha +1}} \\
= & \frac{\alpha \lambda}{2} \int_0^\infty (2 - e^{e^{i\frac{\pi}{2}}y} -e^{e^{-i\frac{\pi}{2}}y}) \frac{dy}{y^{\alpha +1}}\\
= & \frac{\alpha \lambda}{2} \int_0^\infty (1 - e^{e^{i\frac{\pi}{2}}y} ) \frac{dy}{y^{\alpha +1}} + \int_0^\infty (1 - e^{e^{-i\frac{\pi}{2}}y}) \frac{dy}{y^{\alpha +1}}\\
= & \frac{\lambda \Gamma(1-\alpha)}{2} \left( e^{i\frac{\pi \alpha}{2}} +  e^{-i\frac{\pi \alpha}{2}} \right)\\
= & \lambda \Gamma(1-\alpha) \cos \frac{\pi \alpha}{2}>0, \quad \alpha \in (0,1)
\end{align*} 
where we have used \eqref{fur-symb-subordinator}.

\section{Proof of Theorem \ref{second-the}} 
We have that
\begin{align*}
\mathbb{E} e^{i \boldsymbol{\xi} \cdot \mathbf{A}(t)} = & \mathbb{E} \prod_{j=0}^{N(t)} \mathbb{E} e^{i  \epsilon_j \boldsymbol{\xi} \cdot \mathbf{Y}_j } =  \mathbb{E} \Bigg( \mathbb{E} e^{i \epsilon \boldsymbol{\xi} \cdot \mathbf{Y}} \Bigg)^{N(t)} 
\end{align*}
where we used the fact that $\mathbf{Y}_j \sim \mathbf{Y}$ for all $j$. Therefore, we obtain that
\begin{align}
\mathbb{E} e^{i \boldsymbol{\xi} \cdot \mathbf{A}(t)} = & \exp\left( \lambda t (\mathbb{E} e^{i \epsilon \boldsymbol{\xi} \cdot \mathbf{Y}} - 1) \right)\notag \\
= & \exp \left( \frac{\lambda t}{2} (\mathbb{E} e^{i\boldsymbol{\xi} \cdot \mathbf{Y}} + \mathbb{E} e^{-i\boldsymbol{\xi} \cdot \mathbf{Y}} - 2) \right)\label{fur-A-second}
\end{align}
and $\mathbf{Y} \sim \nu_{\mathbf{Y}}$, see formula \eqref{nu-Y-bold}. The Fourier symbol corresponding to the characteristic function \eqref{fur-A-second} is given by
\begin{align*}
-\partial_t\, \mathbb{E} e^{i \boldsymbol{\xi} \cdot \mathbf{A}(t)} \Bigg|_{t=0} = & -\frac{\lambda}{2} \int_{\mathbb{R}^d} (e^{i \boldsymbol{\xi} \cdot \mathbf{y}} + e^{-i \boldsymbol{\xi} \cdot \mathbf{y}} - 2) \nu_{\mathbf{Y}}(\mathbf{dy})\\
= & -\lambda \int_{\mathbb{R}^d} \Big( \cos \boldsymbol{\xi} \cdot \mathbf{y} - 1 \Big) \frac{\gamma^\alpha \, \Gamma(\alpha + \frac{d}{2})}{\pi^\frac{d}{2}\, \Gamma(\alpha)\, (|\mathbf{y}|^2 + \gamma)^{\alpha + \frac{d}{2}}} \mathbf{dy}.
\end{align*}
and therefore, for the process \eqref{A-second-thm}, we get that
\begin{align*}
-\partial_t\, \mathbb{E} e^{i \boldsymbol{\xi} \cdot \mathbf{A}(t/\gamma^\alpha)} \Bigg|_{t=0} =  -\lambda  \int_{\mathbb{R}^d} \Big( \cos \boldsymbol{\xi} \cdot \mathbf{y} - 1 \Big) \frac{ \Gamma(\alpha + \frac{d}{2})}{\pi^\frac{d}{2}\, \Gamma(\alpha)\, (|\mathbf{y}|^2 + \gamma)^{\alpha + \frac{d}{2}}} \mathbf{dy} = \Phi_\gamma(\boldsymbol{\xi}).
\end{align*}
The limit for $\gamma \to 0$ leads to the Fourier symbol
\begin{equation}
\lim_{\gamma \to 0} \Phi_\gamma(\boldsymbol{\xi}) =  -\lambda \frac{ \Gamma(\alpha + \frac{d}{2})}{\pi^\frac{d}{2}\, \Gamma(\alpha)} \int_{\mathbb{R}^d} \Big( \cos \boldsymbol{\xi} \cdot \mathbf{y} - 1 \Big)  \frac{\mathbf{dy}}{|\mathbf{y}|^{2\alpha + d}} = \lambda \frac{ \Gamma(\alpha + \frac{d}{2})}{\pi^\frac{d}{2}\, \Gamma(\alpha)} C\, |\boldsymbol{\xi}|^{2\alpha}  \label{eval-int-1}
\end{equation}
where
\begin{equation}
C =  \int_{\mathbb{R}^d} \frac{1- \cos y_1}{|\mathbf{y}|^{2\alpha + d}} \mathbf{dy}. \label{eval-int-2}
\end{equation}
The interested reader can find  in \cite{NPV2011} a detailed computation of the integrals in \eqref{eval-int-1} and \eqref{eval-int-2}. Finally, we observe that
\begin{align*}
-\partial_t\, \mathbb{E} e^{i \boldsymbol{\xi} \cdot \mathbf{A}(t/\gamma^\alpha)} \Bigg|_{t=0} = & -\frac{\lambda}{2 \gamma^\alpha} \int_{\mathbb{R}^d} (e^{i \boldsymbol{\xi} \cdot \mathbf{y}} + e^{-i \boldsymbol{\xi} \cdot \mathbf{y}} - 2) \nu_{\mathbf{Y}}(\mathbf{dy}) = \Phi_\gamma(\boldsymbol{\xi})
\end{align*}
converges, for $\gamma \to 0$, to the Fourier symbol
\begin{align*}
\Phi(\boldsymbol{\xi}) = -\frac{\lambda}{2} \frac{ \Gamma(\alpha + \frac{d}{2})}{\pi^\frac{d}{2}\, \Gamma(\alpha)} \int_{\mathbb{R}^d} (e^{i \boldsymbol{\xi} \cdot \mathbf{y}} + e^{-i \boldsymbol{\xi} \cdot \mathbf{y}} - 2)  \frac{\mathbf{dy}}{|\mathbf{y}|^{2\alpha + d}} .
\end{align*}
By applying formula \eqref{inv-fourier-symbol}, we get
\begin{align*}
\mathcal{A}f(\mathbf{x}) = & \lambda \frac{ \Gamma(\alpha + \frac{d}{2})}{\pi^\frac{d}{2}\, \Gamma(\alpha)} C \frac{C_d(\alpha)}{2}\int_{\mathbb{R}^d} \Big( f(\mathbf{x} + \mathbf{y}) + f(\mathbf{x} - \mathbf{y}) - 2f(\mathbf{x}) \Big) \frac{ \mathbf{dy}}{|\mathbf{y}|^{2\alpha + d}} \\
= & - \lambda \frac{ \Gamma(\alpha + \frac{d}{2})}{\pi^\frac{d}{2}\, \Gamma(\alpha)} C (-\triangle)^\alpha f(\mathbf{x}) 
\end{align*}
where
\begin{equation*}
C_d(\alpha) = \left( \int_{\mathbb{R}^d} \frac{1- \cos y_1}{|\mathbf{y}|^{2\alpha + d}} \mathbf{dy} \right)^{-1}
\end{equation*}
which is the generator of the isotropic stable process  $\boldsymbol{\mathfrak{S}}^{2\alpha}(t^*)$ with
\begin{equation*}
t^* = \lambda \frac{ \Gamma(\alpha + \frac{d}{2})}{\pi^\frac{d}{2}\, \Gamma(\alpha)} C t, \quad t > 0.
\end{equation*}

\end{document}